\documentclass[a4paper,11pt]{article}
\usepackage{amsmath,amscd,amsfonts,amssymb,epsf,latexsym}

\makeatletter

\long\def\@makefnt#1{\parindent 1em\noindent
            \hb@xt@1.8em{\hss\@textsuperscript{}}#1}
\long\def\@ftntext#1{\insert\footins{%
    \reset@font\footnotesize
    \interlinepenalty\interfootnotelinepenalty
    \splittopskip\footnotesep
    \splitmaxdepth \dp\strutbox \floatingpenalty \@MM
    \hsize\columnwidth \@parboxrestore
    \color@begingroup
      \@makefnt{%
        \rule\z@\footnotesep\ignorespaces#1\@finalstrut\strutbox}%
    \color@endgroup}}%
\def\subjclass#1{%
  \@ftntext{2000 {\itshape Mathematics Subject Classification.}\enspace #1.}}
\def\keywords#1{%
  \@ftntext{{\itshape Key words and phrases.}\enspace #1.}}
\makeatother

\def\A{{\mathbb A}}
\def\B{{\mathbb B}}

\def\C{{\mathbb C}}
\def\D{{\mathbb D}}

\def\X{{\mathbb X}}
\def\Y{{\mathbb Y}}
\def\E{{\mathbb E}}

\def\AB{ {\mathbb A}\moins {\mathbb B}}

\def\moins{\raise 1pt\hbox{{$\scriptstyle -$}}}
\def\plus{\raise 1pt\hbox{{$\scriptstyle +$}} }

\def\phi{\varphi}

\newtheorem{theorem}{Theorem}
\newtheorem{proposition}[theorem]{Proposition}
\newtheorem{lemma}[theorem]{Lemma}
\newtheorem{corollary}[theorem]{Corollary}

\newtheorem{definition}[theorem]{Definition}

\newtheorem{convention}[theorem]{Convention}

\def\proof{\noindent{\bf Proof.\ }}

\def\qed{~\hbox{$\Box$}}

\def\Aut{\mathop{\rm Aut}}
\def\Diff{\mathop{\rm Diff}}
\def\codim{\mathop{\rm codim}}

\begin{document}

\title{\bf Thom polynomials and Schur functions:\\
the singularities $I_{2,2}(-)$}

\author{Piotr Pragacz\thanks{Research supported by a KBN grant,
by T\"UB\.ITAK (during the stay at the METU
in Ankara), and by the Humboldt Stiftung (during the stay at
the MPIM in Bonn).}\\
\small Institute of Mathematics of Polish Academy of Sciences\\
\small \'Sniadeckich 8, 00-956 Warszawa, Poland\\
\small P.Pragacz@impan.gov.pl}



\date{(29.03.2006; revised 15.03.2007)}

\maketitle

\centerline{\it To the memory of Professor Stanis\l aw Balcerzyk
(1932-2005)}

\begin{abstract}
We give the Thom polynomials for the singularities $I_{2,2}$
associated with maps $({\bf C}^{\bullet},0)
\to ({\bf C}^{\bullet+k},0)$ with parameter $k\ge 0$.
Our computations combine the characterization of Thom polynomials
via the ``method of restriction equations'' of Rimanyi et al. with
the techniques of Schur functions.
\end{abstract}

\section{Introduction}\label{intro}

The global behavior of singularities is governed by their {\it Thom
polynomials} (cf. \cite{T}, \cite{Kl}, \cite{AVGL}, \cite{Ka},
\cite{Rim2}). Knowing the Thom polynomial of a singularity $\eta$,
denoted ${\cal T}^{\eta}$, one can compute the cohomology class
represented by the $\eta$-points of a map.
We do not attempt here to survey all activities related
to computations of Thom polynomials, which are difficult tasks in
general.

In the present paper, following a series of papers by Rimanyi et al.
\cite{RS}, \cite{Rim2}, \cite{FR}, \cite{BFR},
we study the Thom polynomials for the singularities
$I_{2,2}$ of the maps $({\bf C}^{\bullet},0) \to ({\bf C}^{\bullet+k},0)$
with parameter $k\ge 0$.

The way of obtaining the thought Thom polynomial is through the solution
of a system of linear equations, which is fine when we want to find
one concrete Thom polynomial, say, for a fixed $k$. However, if we want
to find the Thom polynomials for a series of singularities, associated
with maps $({\bf C}^{\bullet},0) \to ({\bf C}^{\bullet+k},0)$
with $k$ as a parameter, we have to solve {\it simultaneously} a countable
family of systems of linear equations. As stated by Rimanyi
in \cite{Rim2}, p. 512~:

\smallskip

{\it ``However, another challenge is to find Thom polynomials containing
$k$ as a parameter.''}

\smallskip

We do it here for the restriction equations for the singularities $I_{2,2}$
(any $k$) with the help of Schur functions. It appears that a use of
these functions puts a transparent structure on computations of Thom
polynomials. In particular, we get in this way some
recursive formulas (cf., e.g., Lemma \ref{Lr}) that are not so easy to find
using other bases (e.g. the Chern monomial basis that was used in
\cite{Rim2} and in the references quoted there). In fact, various
recursions play a prominent role in the present paper -- apart from
Lemma \ref{Lr}, see Eq.~(\ref{rd}). Let us note that in a recent paper
\cite{FK}, Feher and Komuves compute Thom polynomials for some second
order Thom-Boardman singularities also using Schur functions, and obtain
similar recursions for the coefficients.

Another feature of using the Schur function expansions for Thom
polynomials is that in all known to us cases, all the coefficients
are nonnegative. In fact, we state the following ``positivity conjecture'':

\bigskip

\noindent
{\bf Conjecture:} {\it The coefficients of the Schur function expansion
of a Thom polynomial are nonnegative.}\footnote{Note added in May 2006:
this conjecture (formulated also by Feher and Komuves \cite{FK})
has been recently proved by Weber and the author in \cite{PW}.}

\bigskip

\noindent
To be more precise, we use here (the specializations of)
{\it supersymmetric} Schur functions, also called
``Schur functions in difference of alphabets'' together with their
three basic properties: {\it vanishing}, {\it cancellation} and {\it
factorization}, (cf. \cite{Ste}, \cite{BR}, \cite{LS}, \cite{P2},
\cite {PT}, \cite{M}, \cite{FP}, and \cite{L}). These functions contain
resultants among themselves.
Their geometric role was illuminated, e.g., in the study of
{\it ideals of polynomials supported on degeneracy loci} of \cite{P}, i.e.
${\cal P}$-ideals of singularities $\Sigma^i$ in the terminology of
of the present paper (cf. the end of Section 2 and Theorem 11).

\smallskip

The main goal of this paper is to give the Thom polynomials
for the singularities $I_{2,2}$ (in Mather's notation) associated with
maps $({\bf C}^{\bullet},0) \to ({\bf C}^{\bullet+k},0)$ with parameter
$k\ge 0$. We do it via establishing the Schur function expansions
for these Thom polynomials. We prove first in Lemma \ref{L3} that the
partitions appearing nontrivially have not more than $3$ parts.
Then, in Lemma \ref{Lr}, we establish a recursive relation
for Thom polynomials associated with the successive values of the
parameter $k$. This reduces our calculation to compute the (sub)sum
indexed by partitions with precisely $2$ parts. This is essentially
done in Proposition \ref{Po} (see also Propositions \ref{Pr}, \ref{P<>},
\ref{P<>2}).

Our main result (Theorem \ref{TI22}), combined with Propositions
\ref{P<>} and \ref{P<>2},
gives an explicit presentation of the Thom polynomial for the
singularities $I_{2,2}$ with parameter $k\ge 0$ as a {\bf Z}-combination
of Schur functions. We give {\it closed algebraic} expressions
for the coefficients of these expansions. It turns out that
these coefficients
are the same as the coefficients of the Schur function expansions
of the Segre classes of the second symmetric power of a rank $2$
vector bundle, computed in \cite{S}, \cite{P}, \cite{LLT},
and \cite{P223}.

Our main result offers a generalization (to any $k\ge 0$)
of the formulas obtained previously by Porteous \cite{Po} and Rimanyi
\cite{Rim2} for $k=0$ and $k=1$, respectively.

\smallskip

In our calculations, we use extensively the functorial $\lambda$-ring
approach to symmetric functions developed mainly in
Lascoux's book \cite{L}.

\medskip

Main results of the present paper were announced in \cite{P23}.

\smallskip

The forthcoming author's articles \cite{P3} and \cite{P4} will be devoted
to study the Schur function expansions of Thom polynomials for the Morin's
singularities $A_i$. Inspired by the present paper, \cite{P23}, \cite{P3},
and \cite{P4}, Ozer Ozturk \cite{O} computed the Thom polynomials
for $A_4$ and $k=2,3$.

\section{Recollections on Thom polynomials}

Our main reference for this section is \cite{Rim2}.
We start with recalling what we shall mean by a ``singularity''.
Let $k\ge 0$ be a fixed integer. By a {\it singularity}
we shall mean an equivalence class of stable germs $({\bf C}^{\bullet},0)
\to ({\bf C}^{\bullet+k},0)$, where $\bullet\in {\bf N}$, under the
equivalence
generated by right-left equivalence (i.e. analytic reparametrizations
of the source and target) and suspension (by suspension of a germ $\kappa$
we mean its trivial unfolding: $(x,v) \mapsto (\kappa(x),v)$).

We recall\footnote{This statement is usually called the Thom-Damon
theorem \cite{T}, \cite{D}.}
that the {\it Thom polynomial} ${\cal T}^{\eta}$ of a singularity
$\eta$ is a polynomial in the formal variables $c_1, c_2, \ldots$ that
after the substitution
\begin{equation}
c_i=c_i(f^*TY-TX)=[c(f^*TY)/c(TX)]_i \,,
\end{equation}
for a general map $f:X \to Y$ between complex analytic manifolds,
evaluates the Poincar\'e dual of $[V^{\eta}(f)]$, where $V^{\eta}(f)$
is the cycle carried by the closure of the set
\begin{equation}
\{x\in X : \hbox{the singularity of} \ f \ \hbox{at} \ x \
\hbox{is} \ \eta \}\,.
\end{equation}
By {\it codimension of a singularity} $\eta$, $\codim(\eta)$,
we shall mean $\codim_X(V^{\eta}(f))$ for such an $f$. The concept of
the polynomial ${\cal T}^{\eta}$ comes from Thom's fundamental paper
\cite{T}.
For a detailed discussion of the {\it existence} of Thom polynomials,
see, e.g., \cite{AVGL}. Thom polynomials associated with group actions
were studied by Kazarian in \cite{Ka}, \cite{Ka2}.

In fact, the above is the case with singularities without moduli
\cite{Ka}, \cite{FR}, \cite{BFR}. The singularities $I_{2,2}$, studied
in the present paper, have this property for $k\ge 0$. Indeed, the moduli
of singularities start at codimension greater than $6k+8$ ({\it loc.cit.})
whereas $\codim(I_{2,2})=3k+4$.

According to Mather's classification, singularities are in one-to-one
correspondence with finite dimensional ${\bf C}$-algebras.
We shall use the following notation:

\medskip

-- $A_i$ \ (of Thom-Boardman type $\Sigma^{1_i}$) will stand for the
stable germs with local algebra ${\bf C}[[x]]/(x^{i+1})$, $i\ge 0$;

\medskip

-- $I_{2,2}$ \ (of Thom-Boardman type $\Sigma^2$) for stable germs with
local algebra ${\bf C}[[x,y]]/(xy, x^2+y^2)$;

\medskip

-- $III_{2,2}$ \ (of Thom-Boardman type $\Sigma^2$) for stable germs
with local algebra ${\bf C}[[x,y]]/(xy, x^2, y^2)$\,,
here $k\ge 1$.

\medskip

In the present article, the computations of Thom polynomials shall use
the method which stems from a sequence of papers
by Rimanyi et al. \cite{RS}, \cite{Rim2}, \cite {FR}, \cite{BFR}.
We sketch briefly this approach, refering the interested reader for more
details to these papers, the main references being the last three
mentioned items.

Let $k\ge 0$ be a fixed integer, and let $\eta: ({\bf C}^{\bullet},0) \to
({\bf C}^{\bullet+k},0)$ be a singularity with a prototype
$\kappa: ({\bf C}^n,0) \to ({\bf C}^{n+k},0)$. The {\it maximal compact
subgroup of the right-left symmetry group}
\begin{equation}
\Aut \kappa = \{(\phi,\psi) \in \Diff({\bf C}^n,0) \times
\Diff({\bf C}^{n+k},0) : \psi \circ \kappa \circ \phi^{-1} = \kappa \}
\end{equation}
of $\kappa$ will be denoted by $G_\eta$.
Even if $\Aut \kappa$ is much too large to be a finite dimensional
Lie group, the concept of its maximal compact subgroup (up to conjugacy)
can be defined in a sensible way (cf. \cite{J} and \cite{W}).
In fact, $G_\eta$ can be chosen so that the images of its projections
to the factors $\Diff({\bf C}^n,0)$ and $\Diff({\bf C}^{n+k},0)$ are
linear. Its representations via the projections on the source ${\bf C}^n$
and the target ${\bf C}^{n+k}$ will be denoted by $\lambda_1(\eta)$ and
$\lambda_2(\eta)$.
The vector bundles associated with the universal principal
$G_\eta$-bundle $EG_\eta \to BG_\eta$ using the representations
$\lambda_1(\eta)$ and $\lambda_2(\eta)$ will be called
$E_{\eta}'$ and $E_{\eta}$. The {\it total Chern
class of the singularity} $\eta$ is defined
in $H^{*}(BG_\eta,{\bf Z})$ by
\begin{equation}
c(\eta):=\frac{c(E_{\eta})}{c(E_{\eta}')}\,.
\end{equation}
The {\it Euler class} of $\eta$ is defined in
$H^{2\codim(\eta)}(BG_\eta,{\bf Z})$ by
\begin{equation}
e(\eta):=e(E_{\eta}')\,.
\end{equation}

In the following theorem, we collect information from \cite{Rim2},
Theorem 2.4 and \cite{FR}, Theorem 3.5, needed for the calculations
in the present paper.

\begin{theorem}\label{TEq} \ Suppose, for a singularity $\eta$, that
the Euler classes of all singularities of smaller codimension than
$\codim(\eta)$, are not zero-divisors \footnote{This is the so-called
``Euler condition'' ({\it loc.cit.}). The Euler condition holds true for the
singularities $I_{2,2}$ for any $k\ge 0$.}. Then we have

\noindent
(i) \ if \ $\xi\ne \eta$ \ and \ $\codim(\xi)\le \codim(\eta)$, then
\ ${\cal T}^{\eta}(c(\xi))=0$;

\noindent
(ii) \ ${\cal T}^{\eta}(c(\eta))=e(\eta)$.

\noindent
This system of equations (taken for all such $\xi$'s) determines the
Thom polynomial ${\cal T}^{\eta}$ in a unique way.
\end{theorem}
To use this method of determining the Thom polynomials for
singularities, one needs their classification, see, e.g., \cite{dPW}.

\smallskip

In Section 4, we shall use these equations to compute
Thom polynomials. Sometimes, it is convenient not to work with
the whole maximal compact subgroup $G_\eta$ but with its suitable
subgroup;
this subgroup should be, however, as ``close'' to $G_\eta$ as possible
(cf. \cite{Rim2}, p. 502). We shall denote this subgroup by the
same symbol $G_\eta$.
We recall the following recipe for computing maximal compact subgroups  
from \cite{Rim2} pp. 505--507. Let $\eta$ be a singularity whose
prototype is $\kappa: ({\bf C}^n,0)\to ({\bf C}^{n+k},0)$. The germ
$\kappa$ is the miniversal unfolding of another germ $\beta:
({\bf C}^m,0)\to ({\bf C}^{m+k},0)$ with $d\beta=0$. The group $G_\eta$
is a subgroup of the maximal compact subgroup of
the algebraic automorphism group of
the local algebra $Q_\eta$
of $\eta$ times the unitary group $U(k\moins d)$, where $d$
is the difference between the minimal number of relations and the number
of generators of $Q_\eta$.
With $\beta$ well chosen, $G_\eta$ acts as right-left symmetry group
on $\beta$ with representations $\mu_1$ and $\mu_2$. The representations
$\lambda_1$ and $\lambda_2$ are
\begin{equation}
\lambda_1=\mu_1\oplus \mu_V \ \ \hbox{and} \ \
\lambda_2=\mu_2\oplus \mu_V\,,
\end{equation}
where $\mu_V$ is the representation of $G_\eta$ on the unfolding space
$V={\bf C}^{n-m}$ given, for $\alpha \in V$ and $(\phi,\psi)\in G_\eta$,
by
\begin{equation}
(\phi,\psi) \ \alpha = \psi \circ \alpha \circ \phi^{-1}\,.
\end{equation}
For example, for the singularity of type $A_i$: $({\bf C}^{\bullet},0) \to
({\bf C}^{\bullet+k},0)$, we have $G_{A_i}=U(1)\times U(k)$ with
\begin{equation}
\mu_1=\rho_1, \ \ \mu_2=\rho_1^{i+1}\oplus \rho_k, \ \
\mu_V=\oplus_{j=2}^i \ \rho_1^j \oplus \oplus_{j=1}^i (\rho_k \otimes
\rho_1^{-1})\,,
\end{equation}
where $\rho_j$ denotes the standard representation of
the unitary group $U(j)$. Hence, we obtain assertion (i) of the following

\begin{proposition}\label{Pce}
(i) \ Let $\eta=A_i$. For any $k$, writing $x$ and $y_1$,\ldots, $y_k$
for the Chern roots of the universal bundles on $BU(1)$ and $BU(k)$,
\begin{equation}
c(A_i)=\frac{1+(i+1)x}{1+x}\prod_{j=1}^k(1+y_j)\,.
\end{equation}

\smallskip

\noindent
(ii) \ Let $\eta=I_{2,2}$. Denote by $H$ the extension
of $U(1)\times U(1)$ by ${\bf Z}/2{\bf Z}$ (``the group generated
by multiplication on the coordinates and their exchange'').
For $k\ge 0$, $G_\eta= H \times U(k)$. Hence, for the purpose of our
computations, we can use $G_\eta=U(1)\times U(1) \times U(k)$.
Writing $x_1, x_2$ and $y_1,\ldots, y_k$ for the Chern roots
of the universal bundles
on two copies of $BU(1)$ and on $BU(k)$, we have
\begin{equation}
c(I_{2,2})=\frac{(1+2x_1)(1+2x_2)}{(1+x_1)(1+x_2)}\prod_{j=1}^k(1+y_j)\,,
\end{equation}
\begin{equation}
e(I_{2,2})=x_1x_2(x_1-2x_2)(x_2-2x_1)\prod_{j=1}^k (x_1-y_j)(x_2-y_j)
(x_1+x_2-y_j)\,.
\end{equation}

\smallskip

\noindent
(iii) \ Let $\eta=III_{2,2}$. For $k\ge 1$, $G_\eta=U(2)\times
U(k\moins 1)$,
and
writing $x_1, x_2$ and $y_1,\ldots,y_{k-1}$ for the Chern roots of the
universal bundles on $BU(2)$ and $BU(k\moins 1)$,
\begin{equation}\label{cIII}
c(III_{2,2})=\frac{(1\plus 2x_1)(1\plus 2x_2)(1\plus x_1\plus x_2)}
{(1\plus x_1)(1\plus x_2)}\prod_{j=1}^{k-1}(1+y_j)\,.
\end{equation}
\end{proposition}
(Assertions (ii) and (iii) are obtained, in a standard way, following
the instructions of \cite{Rim2}, Sect.~4. As for assertion (ii), compare
\cite[pp.~506--507]{Rim2} whereas assertion (iii) stems from
\cite[p.~65]{BFR}.)

\smallskip

Let $\eta$ be a singularity. As it was illuminated in the author's
paper \cite{P}, in the case of the singularities $\eta=\Sigma^i$,
it is natural and useful to consider a certain (homogeneous) ideal
in the polynomial ring $R={\bf Z}[c_1, c_2, \ldots]$ whose component
of minimal degree is generated by ${\cal T}^{\eta}$. Namely,
we denote by \ ${\cal P}^{\eta}$ \
the ideal of polynomials in $R$ which -- after the substitution (1) --
are supported on
$V^{\eta}(f)$, where $f: X\to Y$ is a general map
between complex analytic manifolds. (The notion of a ``polynomial
supported on a subscheme'' can be found in \cite{FP}, Appendix A.)
Keeping track of \cite{P}, we shall call ${\cal P}^{\eta}$
the {\it ${\cal P}$--ideal of the singularity $\eta$}.
For example, the ${\cal P}$--ideal of the singularity
$$
\Sigma^i: ({\bf C}^m,0)\to ({\bf C}^n,0)
$$
is
$$
{\cal P}^{\Sigma^i}={\cal P}_{m \moins i}\,,
$$
where on the RHS we have the ideal studied extensively in \cite{P}
(cf. also \cite{P1}, \cite{P2}).
We shall use this ideal in the proof of Theorem \ref{Tp}.

\smallskip

In the present paper, it will be more handy to use, instead of $k$,
a ``shifted'' parameter
\begin{equation}
r:=k+1\,.
\end{equation}
Sometimes, we shall write $\eta(r)$ for the singularity
$\eta: ({\bf C}^{\bullet},0) \to ({\bf C}^{\bullet + r-1},0)$,
and denote the Thom polynomial of $\eta(r)$ by
${\cal T}^{\eta}_r$ -- to emphasize the dependence of both items on $r$.

We have
$$
c_i(f^*TY-TX)=S_i(TX^*-f^*(TY^*))\,,
$$
where $S_i$ means the Segre class. We shall follow of the notation on the
RHS and use, more generally, Schur functions $S_{(i_1,i_2,\ldots,i_h)}$,
indexed by partitions, cf. the next section.

\section{Recollections on Schur functions}

In this section, we collect needed notions related to symmetric functions.
We adopt the functorial point of view
of \cite{L} for what concerns symmetric functions. Namely, given a
commutative ring, we treat symmetric functions as operators acting
on the ring. (Here, these commutative rings are mostly ${\bf Z}$-algebras
generated by the Chern roots of the vector bundles from Proposition
\ref{Pce}.)

\begin{definition}\label{alph}
By an {\it alphabet} $\A$, we understand a (finite)
multi-set of elements in a commutative ring.
\end{definition}
For $m\in {\bf N}$, by ``an alphabet $\A_m$'' we shall mean an alphabet
$\A=(a_1,\ldots,a_m)$ \ (of cardinality $m$); ditto for
$\B_n=(b_1,\ldots,b_n)$, $\Y_k=(y_1,\ldots,y_k)$, and $\X_2=(x_1,x_2)$.

\begin{definition}\label{cf}
Given two alphabets $\A$, $\B$, the {\it complete functions} $S_i(\AB)$
are defined by the generating series (with $z$ an extra variable):
\begin{equation}
\sum S_i(\AB) z^i =\prod_{b\in \B} (1\moins bz)/\prod_{a\in \A}
(1\moins az)\,.
\end{equation}
\end{definition}
So $S_i(\A-\B)$ interpolates between $S_i(\A)$ -- the complete
homogeneous symmetric function of degree $i$ in $\A$
and $S_i(-\B)$ -- the $i$th elementary function in $\B$ times $(-1)^i$.
The notation $\A -\B$ is compatible with the multiplication
of series
\begin{equation}
\sum S_i(\A - \B)z^i \cdot \sum S_j(\A' - \B')z^j =
\sum S_i\bigl((\A+\A') - (\B+\B')\bigr)z^i\,,
\label{}
\end{equation}
the sum $\A + \A'$ denoting the union of two alphabets $\A$ and $\A'$.

\begin{convention} \rm We shall often identify an alphabet
$\A=\{a_1,\ldots,a_m\}$ with the sum $a_1+\cdots +a_m$ and perform usual
algebraic operations on such elements. For example, $\A b$ will
denote the alphabet $(a_1b,\ldots,a_mb)$.
We will give priority to the
algebraic notation over the set-theoretic one. In fact, in the following,
we shall use mostly alphabets of variables.
\end{convention}

We have $(\A+\C) - (\B+\C) = \A-\B$, and this corresponds to
simplification of the common factor for the rational series:
\begin{equation}\label{Canc}
\sum S_i((\A + \C) - (\B + \C))z^i = \sum S_i(\A-\B)z^i\,.
\end{equation}

\begin{definition}
By a partition $I=(i_1,i_2,\ldots,i_h)$ we mean a weakly increasing
sequence $0\le i_1\le i_2\le \ldots \le i_h$ of natural numbers.
\end{definition}
In the following, we shall identify partitions with their {\it Young
diagrams}, as is customary.

\begin{definition}
Given a partition $I$ and two alphabets $\A$ and $\B$,
the {\it Schur function} $S_I(\A \moins \B)$ is defined by
the following determinant:
\begin{equation}\label{schur}
S_I(\A \moins \B):= \Bigl|
     S_{i_p+p-q}(\A \moins \B) \Bigr|_{1\leq p,q\le h}  \ .
\end{equation}
\end{definition}
These functions are often called {\it supersymmetric Schur functions}
or {\it Schur functions in difference of alphabets}. Their properties
were studied, among others, in \cite{Ste}, \cite{BR}, \cite{LS}, \cite{P2},
\cite {PT}, \cite{M}, \cite{FP}, and \cite{L}. From the last item,
we borrow a use of increasing ``French'' partitions and the determinant of
the form (\ref{schur}) evaluating a Schur function.
We shall use the simplified notation $i_1i_2\cdots i_h$ or
$i_1,i_2,\ldots , i_h$ for a partition $(i_1,i_2,\ldots,i_h)$
(the latter one if $i_h\ge 10$). The rectangle partition $(i,i,\ldots,i)$
($h$ times) will be denoted $(i^h)$.

For example,
$$
S_{33344}(\AB)
=\begin{vmatrix}
S_3 \ & \ S_4 \ & S_5 \ & \ S_7 \ &  \ S_8 \\
S_2 & S_3 & S_4 & S_6 & S_7 \\
S_1 & S_2 & S_3 & S_5 & S_6 \\
1 & S_1 & S_2 & S_4 & S_5 \\
0 & 1 & S_1 & S_3 & S_4
\end{vmatrix}\,,
$$
where $S_i$ means $S_i(\AB)$.

We shall now give some properties of Schur functions. The details can
be found in the just quoted references.
By Eq.~(\ref{Canc}), we get the following {\it cancellation property}:
\begin{equation}
S_I((\A + \C) - (\B + \C))=S_I(\A-\B)\,.
\end{equation}

We record the following property justifying the notational remark
from the end of Section 2; for a partition $I$,
\begin{equation}
S_I(\AB)= (-1)^{|I|}S_J(\B \moins \A)=S_J(\B^* \moins \A^*)\,,
\end{equation}
where $J$ is the conjugate partition of $I$ (i.e. the consecutive
rows of $J$ are equal to the corresponding columns of $I$), and
$\A^*$ denotes the alphabet $\{-a_1,-a_2,\ldots \}$.

Fix two positive integers $m$ and $n$.
We shall say that a partition $I=(i_1,i_2,\ldots,i_h)$
{\it is contained in} the $(m,n)$-hook if either $h\le m$, or $h> m$
and $i_{h-m}\le n$.
Pictorially, this means that the Young diagram of $I$ is contained
in the ``tickened" hook

\smallskip

$$
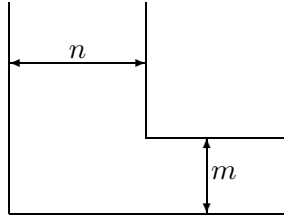
\unitlength=2mm
\begin{picture}(18,14)
\put(0,0){\line(0,1){14}}
\put(0,0){\line(1,0){18}}
\put(9,5){\line(1,0){9}}
\put(9,5){\line(0,1){9}}
\put(4,10){\vector(1,0)5}
\put(4,10){\vector(-1,0)4}
\put(13,2){\vector(0,1)3}
\put(13,2){\vector(0,-1)2}
\put(4.5,10.3){\hbox to0pt{\hss$n$\hss}}
\put(13.3,2.3){\hbox{$m$}}
\end{picture}
$$
We record the following {\it vanishing property}.
Given alphabets $\A$ and $\B$ of cardinalities $m$ and $n$, if
a partition $I$ is not contained in the $(m,n)$-hook,
then
\begin{equation}\label{van}
S_I(\A-\B)=0\,.
\end{equation}
For example,
$$
S_{3569}(\A_2-\B_4)=
S_{3569}(a_1\plus a_2\moins b_1\moins b_2\moins b_3\moins b_4)=0
$$
because $3569$ is not contained in the $(2,4)$-hook.
In fact, we have the following result.
\begin{theorem}\label{Tss}
If $\A_m$ and $\B_n$ are alphabets of variables, then the functions
$S_I(\A_m-\B_n)$, for $I$ running over partitions
contained in the $(m,n)$-hook, are ${\bf Z}$-linearly independent.
\end{theorem}
(They form a ${\bf Z}$-basis of the abelian group of the so-called
``supersymmetric functions''.)

In the present paper, by a {\it symmetric function} we shall mean
a ${\bf Z}$-linear combination of the operators $S_I(-)$.

\begin{definition}
Given two alphabets $\A,\B$, we define their {\it resultant}
\begin{equation}\label{}
R(\A,\B):=\prod_{a\in \A,\, b\in \B}(a\moins b)\,.
\end{equation}
\end{definition}
This terminology is justified by the fact that $R(\A,\B)$ is the classical
resultant of the polynomials $R(x,\A)$ and $R(x,\B)$.
We have
\begin{equation}\label{ER}
R(\A_m,\B_n)= S_{(n^m)}(\AB)=\sum_I S_I(\A) S_{(n^m)/I}(-\B)\,,
\end{equation}
where the sum is over all partitions $I\subset (n^m)$.

When a partition is contained in the $(m,n)$-hook and at the same time it
contains the rectangle $(n^m)$, then we have the following
{\it factorization property}:
for partitions $I=(i_1,\ldots,i_m)$ and $J=(j_1,\ldots, j_h)$,
\begin{equation}\label{Fact}
S_{(j_1,\ldots,j_h,i_1+n,\ldots,i_m+n)}(\A_m-\B_n)
=S_I(\A) \ R(\A,\B) \ S_J(-\B)\,.
\end{equation}

\smallskip

The following convention stems from Lascoux's paper \cite{L1}.

\begin{convention} \rm We may need to specialize a letter to $2$, but this must
not be confused with taking two copies of $1$. To allow one, nevertheless,
specializing a letter to an (integer, or even complex) number $r$ inside
a symmetric function, without introducing intermediate variables,
we write \fbox{$r$} for this specialization. Boxes have to be treated
as single variables. For example, $S_i(2) = {{i+1} \choose 2}$ but
$S_i(\fbox{$2$})=2^i$.
A similar remark applies to ${\bf Z}$-linear combinations of variables.
We have \ $S_2(\X_2)=x_1^2\plus x_1x_2\plus x_2^2$ but
$S_2(\fbox{$x_1\plus x_2$})=x_1^2\plus 2x_1x_2\plus x_2^2$,
$S_{11}(\X_2)=x_1x_2$ but
$S_{11}(\fbox{$x_1\plus x_2$})=0$, $S_{2}(3x)=6x^2$ but
$S_{2}(\fbox{$3x$})=9x^2$ etc.
\end{convention}

\smallskip

For example,
\begin{equation}\label{eR}
\prod_{j=1}^k (x_1-y_j)(x_2-y_j)(x_1+x_2-y_j)=
R(\X_2\plus \fbox{$x_1\plus x_2$},\Y_{k})\,.
\end{equation}

\smallskip

\noindent
This convention will be used in the next section.

\medskip

We end the present section with the following result which is
a consequence of the author's study \cite{P}, \cite{P1}, \cite{P2}
of the ${\cal P}$-ideals of the singularities $\Sigma^i$.

\begin{theorem}\label{Tp}
Suppose that a singularity $\eta$ is of Thom-Boardman type
$\Sigma^i$. Then all summands
in the Schur function expansion of ${\cal T}^{\eta}_r$ are indexed
by partitions containing\footnote{We say that one partition {\it is
contained} in another if this holds for their Young diagrams
(cf. \cite{L}).} the rectangle partition $(r+i-1)^i$.
\end{theorem}
\proof
Since $\eta$ is of Thom-Boardman type $\Sigma^i$, the Thom
polynomial ${\cal T}^{\eta}_r$ belongs to the $\cal{P}$--ideal of
the singularity $\Sigma^i$ with parameter $r$. We also know
by the Thom-Damon theorem (cf. \cite{D}) that ${\cal T}^{\eta}_r$
is a ${\bf Z}$-linear combination of Schur functions in $TX^*-f^*(TY^*)$.
The assertion now follows by combining Theorem 3.4 from \cite{P}
with Lemma 2.5 from \cite{P1} (see also Claim in the proof of
Theorem 5.3(i) in \cite{P2}). Indeed, it follows from the former result
that any ${\bf Z}$-combination of Schur functions indexed by partitions
containing $(r+i-1)^i$ belongs to ${\cal P}^{\Sigma_i(r)}$, whereas
the latter result implies that no nonzero
${\bf Z}$-combination $\sum_I \alpha_I S_I$,
where all $I\not\supset (r+i-1)^i$, belongs to ${\cal P}^{\Sigma_i(r)}$.
\qed

\section{Thom polynomial for $I_{2,2}(r)$}

The codimension of $I_{2,2}(r)$, $r\ge 1$, is $3r+1$.
The Thom polynomial for $I_{2,2}(1)$
is \ $S_{22}=S_{22}(TX^*\moins f^*(TY^*))$ \
(cf. \cite{Po}). In the following, we shall often omit the arguments
of Schur functions.

From now on, we shall assume that $r\ge 2$. The Thom polynomial
for $I_{2,2}(2)$ is (cf. \cite {Rim2})
$$
S_{133}+3S_{34}\,.
$$
By virtue of Proposition \ref{Pce}, the equations from Theorem \ref{TEq}
characterizing the Thom polynomial for $I_{2,2}(r)$ are
\begin{equation}\label{Ia}
P(-\B_{r-1})=P(x-\fbox{$2x$}-\B_{r-1})=P(x-\fbox{$3x$}-\B_{r-1})=0\,,
\end{equation}
and (using Eq. (\ref{eR}))
\begin{equation}\label{Ib}
P(\X_2\moins \fbox{$2x_1$}\moins \fbox{$2x_2$}\moins \B_{r-1})
=x_1x_2(x_1\moins 2x_2)(x_2\moins 2x_1)
\ R(\X_2\plus \fbox{$x_1\plus x_2$},\B_{r-1})\,.
\end{equation}
Here, without loss of generality, we assume that $x$, $x_1$, $x_2$,
and $\B_{r-1}$ are variables.
Moreover, $P(-)$ denotes a symmetric function.
For the remainder of this paper, we set
\begin{equation}
\D:=\fbox{$2x_1$}+\fbox{$2x_2$}+\fbox{$x_1+x_2$}\,.
\end{equation}
Then, additionally, for variables $x_1,x_2$ and an alphabet $\B_{r-2}$,
we have the vanishing imposed by $III_{2,2}$~:
\begin{equation}\label{Ic}
P(\X_2-\D -\B_{r-2})=0\,.
\end{equation}
Indeed, the singularities $\ne I_{2,2}$ with codimension
$\le \codim(I_{2,2})$ are: $A_0$, $A_1$, $A_2$, $III_{2,2}$.

For $r\ge 1$, we set
\begin{equation}
{\cal T}_r :={\cal T}^{I_{2,2}}_r\,.
\end{equation}
Our goal is to give a presentation of ${\cal T}_r$ as a {\bf Z}-linear
combination of Schur functions with explicit algebraic expressions of the
coefficients:
\begin{equation}\label{eT}
{\cal T}_r=\sum_I \alpha_I S_I.
\end{equation}
We shall say that a partition $I$ appears nontrivially in Eq.~(\ref{eT})
if $\alpha_I\ne 0$.

\begin{lemma}\label{L3}
(i) A partition appearing nontrivially in the Schur function expansion
of ${\cal T}_r$ contains the partition $(r+1,r+1)$.

\smallskip

\noindent
(ii) A partition appearing nontrivially in the Schur function expansion
of ${\cal T}_r$ has at most three parts.
\end{lemma}
Proof.
(i) Since the singularity $I_{2,2}$ is of Thom-Boardman type $\Sigma^2$,
this is a particular case of Theorem \ref{Tp}.

\smallskip

\noindent
(ii) We can assume that $r\ge 3$. In addition to information
contained in (i), we shall use Eq.~(\ref{Ic})
$$
{\cal T}_r(\X_2-\D -\B_{r-2})=0\,.
$$
By virtue of (i), we can use factorization property (\ref{Fact})
to all summands of
\begin{equation}
{\cal T}_r(\X_2-\D-\B_{r-2})=\sum_I \alpha_I S_I(\X_2-\D-\B_{r-2})
\end{equation}
(we assume that $\alpha_I\ne 0$). We divide each summand of this last
polynomial by the resultant
$$
R(\X_2, \D+\B_{r-2})\,.
$$
Suppose that the resulting factor of $S_I$ is
\begin{equation}\label{factor}
S_{p,q}(\X_2) \ S_J(-\D-\B_{r-2})\,,
\end{equation}
cf. (\ref{Fact}). Since $|I|=3r+1$, we have
\begin{equation}\label{r-1}
|J|\le r-1 \,.
\end{equation}
Now, let us assume that $I$ has more than 3 parts, that is,
$J$ has 2 or more parts. This assumption (together with the
inequality (\ref{r-1})) implies that
$$
S_J(-\B_{r-2})\ne 0
$$
($\B_{r-2}$ is an alphabet of variables).
Expanding (\ref{factor}), we get among summands the following one of
largest possible degree $|J|$ in $\B_{r-2}$:
\begin{equation}\label{summ}
S_{p,q}(\X_2) \ S_J(-\B_{r-2})\ne 0\,.
\end{equation}
Take in the sum
$$
\sum_I \alpha_I S_{p,q}(\X_2) \ S_J(-\D-\B_{r-2})
$$
the (sub)sum of all the nonzero summands of the form (\ref{factor}) with
the largest possible weight of $J$.
Since Schur polynomials are independent, this (sub)sum is nonzero, and
moreover, it is ${\bf Z}$-linearly independent of other summands both
in the sum indexed by partitions $I$ with $\ge 3$ parts, and in that
indexed by partitions with 2 parts (this last sum does not depend on
$\B_{r-2}$).
Hence, there is no {\bf Z}-linear combination of $S_I$'s involving
nontrivially $I$ with more than 3 parts (and possibly also those with 3
or 2 parts) that satisfies Eq.~(\ref{Ic}). Assertion (ii) has been proved.
\qed

\smallskip

\noindent
(For example, $S_{1144}$ cannot appear in the Schur function
expansion of ${\cal T}_3$ because $S_{1144}(\X_2-\D-\B_1)$ after division
by the resultant contains the summand $S_{11}(-\B_1)=S_2(\B_1)$, which
does not occur in similar expressions for $S_{55}, S_{46}, S_{244},
S_{145}$.)

\begin{corollary}\label{Co}
If $S_{i_1,i_2}$ appears nontrivially in the Schur function expansion
of ${\cal T}_r$, then $i_1=r+1+p$ and $i_2=2r-p$,
where $0\le 2p\le r-1$.
\end{corollary}

\medskip

The following lemma gives a recursive description of ${\cal T}_r$.
Denote by $\Phi$ the linear endomorphism on the ${\bf Z}$-module
spanned by Schur functions indexed by
partitions of length $\le 3$,
that sends a Schur function $S_{i_1,i_2,i_3}$ to $S_{i_1+1,i_2+1,i_3+1}$.
Let $\overline{{\cal T}_r}$ denote the sum of those terms
in the Schur function expansion
of ${\cal T}_r$ which are indexed by partitions of length $\le 2$.
Note that ${\overline{\cal T}_1}={\cal T}_1=S_{22}$.

\begin{lemma}\label{Lr} With this notation, for $r\ge 2$, we have the following
recursive equation:
\begin{equation}\label{lr}
{\cal T}_r={\overline{\cal T}_r}+\Phi({\cal T}_{r-1})\,.
\end{equation}
\end{lemma}
\proof
Write
\begin{equation}\label{dP}
{\cal T}_r=\sum_I \alpha_I S_I=\sum_J \alpha_J S_J + \sum_K \alpha_K S_K\,,
\end{equation}
where $J$ have 2 parts and $K=(k_1,k_2,k_3)$ have 3 parts
(we assume that $\alpha_I\ne 0$). We set
\begin{equation}
Q=\sum_K \alpha_K S_{k_1-1,k_2-1,k_3-1}\,,
\end{equation}
and our goal is to show that $Q={\cal T}_{r-1}$.
Since a partition $I$ appearing nontrivially in the Schur function
expansion
of ${\cal T}_r$ must contain the partition $(r\plus 1,r\plus 1)$, any
partition $(k_1\moins 1,k_2\moins 1,k_3\moins 1)$ above contains
the partition $(r,r)$. Since this last partition
is not contained in the $(1,r-1)$-hook, Eqs.~(\ref{Ia}) with $r$ replaced
by $r-1$ and $P$ by $Q$ are automatically fulfilled by virtue of the
vanishing property (\ref{van}).
Note that Eq.~(\ref{Ic}) is a particular case of Eq.~(\ref{Ib}).
Indeed, specializing $b_{r-1}$ to \fbox{$x_1\plus x_2$} in Eq.~(\ref{Ib}),
we get Eq.~(\ref{Ic}). Therefore it suffices to show that
\begin{equation}
Q(\X_2-\E-\B_{r-2})=
x_1x_2(x_1\moins 2x_2)(x_2\moins 2x_1)
\ R(\X_2\plus \fbox{$x_1\plus x_2$},\B_{r-2})\,,
\end{equation}
where $\E=\fbox{$2x_1$}+\fbox{$2x_2$}$\,.
We apply to each summand
$$
\alpha_K S_{k_1-1,k_2-1,k_3-1}(\X_2-\E-\B_{r-2})
$$
of $Q(\X_2-\E-\B_{r-2})$
the factorization property (\ref{Fact}), and divide it by the resultant
$$
R(\X_2, \E+\B_{r-2})\,.
$$
Suppose that the resulting factor is
\begin{equation}
\alpha_K S_{a,b}(\X_2) \ S_c(-\E-\B_{r-2})\,,
\end{equation}
where $(k_1\moins 1,k_2\moins 1,k_3\moins 1)=(c,r\plus a,r\plus b)$.

Performing the same division of
$$
x_1x_2(x_1\moins 2x_2)(x_2\moins 2x_1)
\ R(\X_2\plus \fbox{$x_1\plus x_2$},\B_{r-2})\,,
$$
we get $R(\fbox{$x_1\plus x_2$},\B_{r-2})$. Thus, the desired
equation $Q={\cal T}_{r-1}$ is equivalent to
\begin{equation}\label{show}
\sum_{a+b+c=r-2} \alpha_K S_{a,b}(\X_2) \ S_c(-\E-\B_{r-2})=
R(\fbox{$x_1\plus x_2$},\B_{r-2})\,.
\end{equation}
To prove Eq.~(\ref{show}), we use Eqs.~(\ref{Ib}) and (\ref{dP})
for ${\cal T}_r$
$$
\sum_I \alpha_I S_I(\X_2\moins \E \moins \B_{r-1})
=x_1x_2(x_1\moins 2x_2)(x_2\moins 2x_1)
\ R(\X_2\plus \fbox{$x_1\plus x_2$},\B_{r-1})\,.
$$
Using again the factorization property (this time w.r.t. the larger
rectangle $(r\plus 1,r\plus 1)$) and dividing both sides of the last
equation by the resultant
$$
R(\X_2, \E+\B_{r-1})\,,
$$
we get the identity
\begin{equation}\label{eqat}
\sum_{p+q+j=r-1} \alpha_I S_{p,q}(\X_2) \ S_j(-\E-\B_{r-1})=
R(\fbox{$x_1\plus x_2$},\B_{r-1}).
\end{equation}
Since
$$
S_j(-\E -\B_{r-1})=S_j(-\E-\B_{r-2})-b_{r-1}S_{j-1}(-\E -\B_{r-2})
$$
and
$$
R(\fbox{$x_1\plus x_2$},\B_{r-1})=(x_1+x_2-b_{r-1})
R(\fbox{$x_1\plus x_2$},\B_{r-2})\,,
$$
taking the coefficients of $(-b_{r-1})$ in both sides of Eq.~(\ref{eqat}),
we get the desired Eq.~(\ref{show}). The lemma has been proved.
\qed

\smallskip

\noindent
(For example, writing ${\cal T}_3=\alpha S_{46}+\beta S_{55}+\gamma S_{244}+\delta
S_{145}$, we get that
$$
\gamma S_1(-\E-B_1)+\delta S_1(\X_2)= R(\fbox{$x_1\plus x_2$}, \B_1)\,,
$$
by taking the coeficients of $(-b_2)$ in both sides of
$$
\alpha S_2(\X_2)+\beta S_{11}(\X_2)+\gamma S_2(-\E-\B_2)+\delta
S_1(-\E-\B_2)S_1(X_2)
= R(\fbox{$x_1\plus x_2$}, \B_2)\,.)
$$

\smallskip

Iterating Eq.~(\ref{lr}) gives
\begin{corollary}\label{Cr} With the above notation, we have
\begin{equation}\label{cr}
{\cal T}_r={\overline{\cal T}_r}+\Phi({\overline{\cal T}_{r-1}})
+\Phi^2({\overline{\cal T}_{r-2}})+\cdots+\Phi^{r-1}
({\overline{\cal T}_1)}\,.
\end{equation}
\end{corollary}

Of course, ${\overline{\cal T}_r}$ is uniquely determined by its value on
$\X_2$. The following result gives this value.

\begin{proposition}\label{Po}
For any $r\ge 1$, we have
\begin{equation}
{\overline{\cal T}_r}(\X_2)=(x_1x_2)^{r+1} \ S_{r-1}(\D)\,.
\end{equation}
\end{proposition}
\proof
We use induction on $r$. For $r=1,2$, the assertion holds true.
Suppose that the assertion is true for ${\overline{\cal T}_i}$
where $i<r$.
Fix a partition $I=(j,r+1+p,r+1+q)$ appearing nontrivially in
the Schur function expansion (\ref{eT}) of ${\cal T}_r$.
Note that $j$ varies from $0$ to $r-1$ because $|I|=3r+1$.
We obtain, by the factorization property (\ref{Fact}),
$$
S_I(\X_2-\D-\B_{r-2})=R \cdot S_j(-\D-\B_{r-2}) \cdot S_{p,q}(\X_2)\,,
$$
where $R=R(\X_2, \D + \B_{r-2})$. Hence, using Eq.~(\ref{cr}),
we see that
\begin{equation}\label{sumj}
{\cal T}_r(\X_2-\D-\B_{r-2}) = R \cdot
\Bigl(\sum_{j=0}^{r-1} S_j(-\D-\B_{r-2})
\ \frac{{\overline{\cal T}_{r-j}}(\X_2)}{(x_1x_2)^{r-j+1}}\Bigr)\,.
\end{equation}
By the induction assumption, for positive $j\le r-1$,
$$
{\overline{\cal T}_{r-j}}(\X_2)=(x_1x_2)^{r-j+1} \ S_{r-1-j}(\D)\,.
$$
Substituting this to (\ref{sumj}), and using the vanishing (\ref{Ic}),
we obtain
\begin{equation}\label{l1}
\sum_{j=1}^{r-1} S_j(-\D-\B_{r-2}) S_{r-1-j}(\D)
+\frac{{\overline{\cal T}_r}(\X_2)}{(x_1x_2)^{r+1}}=0\,.
\end{equation}
But we also have, by a formula for addition of alphabets,
\begin{equation}\label{l2}
\sum_{j=1}^{r-1} S_j(-\D-\B_{r-2}) S_{r-1-j}(\D)+S_{r-1}(\D)
=S_{r-1}(-\B_{r-2})=0\,.
\end{equation}
Combining Eqs.~(\ref{l1}) and (\ref{l2}) gives
$$
{\overline{\cal T}_r}(\X_2)=(x_1x_2)^{r+1} \ S_{r-1}(\D)\,,
$$
that is, the induction assertion.
The proof of the proposition is now complete.
\qed

\medskip

The Schur function expansion of $S_i(\D)$ was described in \cite{P},
\cite{LLT}, and \cite[App. A3]{P223} in the context of the
{\it Segre classes} of the second symmetric power of a rank $2$ vector
bundle. Indeed, $\D$ is the alphabet of the Chern roots of the second
symmetric power of a rank $2$ bundle with the Chern roots $x_1, x_2$.

Denote by $\langle p,q\rangle$ the coefficient of $S_{p,q}:=S_{p,q}(\X_2)$
in $S_{p+q}(\D)$, where $0\le p\le q$. A proof of the next proposition,
due to Lascoux with the help of {\it divided differences}, can be found
in \cite{P223}, pp.~163--166. We give here another proof without divided
differences.

\begin{proposition}\label{Pr}
For $p>0$, we have
\begin{equation}\label{erec}
\langle p,q\rangle =\langle p-1,q\rangle +\langle p,q-1\rangle \,.
\end{equation}
\end{proposition}
\proof
We have
\begin{equation}\label{eSi}
S_i(\D)=
\sum_{h=0}^i S_h(\fbox{$2x_1$}\plus \fbox{$2x_2$})S_{i-h}
(\fbox{$x_1\plus x_2$})
=\sum_{h=0}^i 2^h S_h \cdot (x_1\plus x_2)^{i-h} \,,
\end{equation}
and (cf., e.g., \cite{M} I.4, Ex.3)
\begin{equation}\label{ealpha}
(x_1+x_2)^j=\sum_{a,b\ge 0} {a+b \choose a} \frac{b-a+1}{b+1} S_{a,b}\,,
\end{equation}
where $a+b=j$ and $a\le b$.
Combining Eqs.~(\ref{eSi}), (\ref{ealpha}) with the Pieri formula
(cf., e.g., \cite{L}, \cite{M}), we get for $0\le p \le q$,
\begin{equation}\label{epq}
\langle p,q\rangle =\sum_{h=0}^{p+q} \ 2^h \sum_{h_1,h_2\ge 0}
{p\plus q\moins h \choose p\moins h_1}
\frac{(q\moins h_2) \moins (p\moins h_1)\plus 1}{q\moins h_2\plus 1}\,,
\end{equation}
where $h_1+h_2=h$ and $h_1 \le p \le q\moins h_2$.

We also compute the Schur function expansion of $S_{1,i-1}(\D)$.
Denote by $[p,q]$ the coefficient of $S_{p,q}$ in $S_{1,p+q-1}(\D)$,
$0\le p\le q$. We have the following expansion for $S_{1,i-1}(\D)$:
$$
\aligned
&\sum_{h=1}^i S_{(1,i-1)/(i-h)}(\fbox{$2x_1$}\plus \fbox{$2x_2$})S_h
(\fbox{$x_1\plus x_2$})\\
&=\sum_{h=1}^i 2^h S_h \cdot (x_1\plus x_2)^{i-h}
+\sum_{h=1}^i 2^h S_{1,h-1} \cdot (x_1\plus x_2)^{i-h}\,.
\endaligned
$$
We get from both sums in the last line that for $p>0$ the coefficient
$[p,q]$ is equal {\it twice} the RHS of Eq.~(\ref{epq}), that is,
\begin{equation}\label{e2pq}
[p,q]=2\langle p,q\rangle \,.
\end{equation}

We have, by the Pieri formula,
\begin{equation}
S_{i-1}(\D)\cdot S_1(\D)=S_{i-1}(\D) \cdot 3S_1=S_i(\D)+S_{1,i-1}(\D)\,.
\end{equation}
This equation implies that $S_{p,q}$ appears in $S_i(\D)+S_{1,i-1}(\D)$
with multiplicity $3(\langle p\moins 1,q\rangle +\langle p-1,q\rangle )$
(we use the Pieri formula once again).
The desired Eq.~(\ref{erec}) now follows by virtue of Eq.~(\ref{e2pq}).
\qed

\smallskip

We now pass to some ``closed'' algebraic expressions for the
$\langle p,q\rangle $'s.
We have
\begin{equation}\label{e2^}
\langle 0,q\rangle =S_q(\fbox{$1$}+\fbox{$2$})=1+2+\cdots+2^q=2^{q+1}-1\,.
\end{equation}

The following result was obtained in \cite{S}, \cite{P}, and \cite{LLT}.

\begin{proposition}\label{P<>} For $0\le p\le q$, we have
\begin{equation}\label{Seg}
\langle p,q\rangle ={p+q+1\choose p+1}+{p+q+1\choose p+2}+\cdots
+{p+q+1\choose q+1}\,.
\end{equation}
\end{proposition}
We propose now an alternative expression involving powers of $2$,
which is a natural generalization of the equation
$\langle 0,q \rangle =2^{q+1}\moins 1$,
and which stems directly from Eq.~(\ref{erec}).
Namely, with the convention that ${a \choose 0}=1$
for any $a\in {\bf Z}$, we have
\begin{proposition}\label{P<>2}
For $0\le p\le q$,
\begin{equation}\label{eb}
\langle p, q\rangle =2^{p+q+1}
-\sum_{s=0}^p\Big[{p\plus q\moins 2s\moins 1 \choose
p\moins s}
-{2p\moins 2s\moins 1 \choose p\moins s\plus 1}\Big] 2^{2s} \,.
\end{equation}
\end{proposition}
\proof
The proof uses double induction on $p$ and $q$. We use Eq.~(\ref{erec})
several times:
$$
\aligned
\langle p,q\rangle
&=\langle p\moins 1,q\rangle+\langle p,q\moins 1\rangle \\
&=\langle p\moins 1,q\rangle+\langle p\moins 1,q\moins 1\rangle
+\langle p,q\moins 2\rangle \\
&=\ldots \\
&=\langle p\moins 1,q\rangle+\cdots +\langle p\moins 1,1\rangle
+\langle p,0\rangle \,.
\endaligned
$$
We know the values of all summands in the last row by the induction
assumption (the last summand being equal to $2^{p+1}\moins 1$).
Using several times Eq.~(\ref{e2^}) as well as a well-known equality:
$$
1 + {a\plus 1 \choose a} + {a\plus 2 \choose a}+\cdots +
{2a\moins 2 \choose a}
={2a\moins 1 \choose a\plus 1}\,,
$$
we get the desired induction assertion (\ref{eb}) for $\langle p,q\rangle$.
\qed

\bigskip

Using Proposition \ref{Po}, we shall now give the Schur function expansion
of ${\overline{\cal T}_r}$. Denote by \ $d_{rj}$ \ the coefficient
of $S_{r+j,2r+1-j}$
in ${\overline{\cal T}_r}$ for $r\ge 1$ and $j\ge 1$.
By virtue of Corollary \ref{Co},
$d_{rj}\ne 0$ entails $j \le [(r+1)/2]$ (for example, the only Schur
functions that can appear with nonzero coefficients in ${\overline
{\cal T}_5}$ are
$S_{6,10}$, $S_{79}$, and $S_{88}$), so that we have
\begin{equation}
{\overline{\cal T}_r}=\sum_{j=1}^{[(r+1)/2]} \ d_{rj} \ S_{r+j,2r+1-j}\,.
\end{equation}

\medskip

We have the following link between the $d_{rj}$'s and $\langle p,q\rangle$'s:
suppose that $d_{rj}\ne 0$, then we have
\begin{equation}\label{ed<>}
d_{rj}=\langle j\moins 1,r\moins j \rangle\,.
\end{equation}
We may display the $d_{rj}$'s with the help of the following
``Pascal triangle''-type matrix:

\smallskip

$$
\begin{array}{cccccc}
d_{11} & 0 & 0 & 0 & 0 & \ldots \\
d_{21} & 0 & 0 & 0 & 0 & \ldots \\
d_{31} & d_{32} & 0 & 0 & 0 & \ldots \\
d_{41} & d_{42} & 0 & 0 & 0 & \ldots \\
d_{51} & d_{52} & d_{53} & 0 & 0 & \ldots \\
d_{61} & d_{62} & d_{63} & 0 & 0 & \ldots \\
d_{71} & d_{72} & d_{73} & d_{74} & 0 & \ldots \\
\vdots & \vdots & \vdots & \vdots & \vdots &
\end{array} \ \ \ \ \ \
= \ \ \ \ \ \
\begin{array}{cccccc}
1 & 0 & 0 & 0 & 0 & \ \ldots \\
3 & 0 & 0 & 0 & 0 & \ldots \\
7 & 3 & 0 & 0 & 0 & \ldots \\
15 & 10 & 0 & 0 & 0 & \ldots \\
31 & 25 & 10 & 0 & 0 & \ldots \\
63 & 56 & 35 & 0 & 0 & \ldots \\
127 & 119 & 91 & 35 & 0 & \ldots \\
\vdots & \vdots & \vdots & \vdots & \vdots &
\end{array}
$$
By Proposition \ref{Pr}, if $d_{rj}>0$, then we have
\begin{equation}\label{rd}
d_{rj}= d_{r-1,j-1} + d_{r-1,j}\,.
\end{equation}

\medskip

\noindent
We have the following values of ${\overline{\cal T}_1},
{\overline{\cal T}_2}, \ldots, {\overline {\cal T}_7}$:
$$
\aligned
&S_{22}\\
&3S_{34}\\
&7S_{46}+3S_{55}\\
&15S_{58}+10S_{67}\\
&31S_{6,10}+25S_{79}+10S_{88}\\
&63S_{7,12}+56S_{8,11}+35S_{9,10}\\
&127S_{8,14}+119S_{9,13}+91S_{10,12}+35S_{11,11}\,.
\endaligned
$$

\bigskip

Summing up all our considerations, we get the main result of the
present paper. It gives the desired Thom polynomial in a parametric
form (the parameter being $r$).

\begin{theorem}\label{TI22} For $r\ge 1$, the Thom polynomial
for $I_{2,2}(r)$ is equal to
\begin{equation}
\sum_{i=0}^{r-1} \ \sum_{\{j\ge 1: \ i+2j\le r+1\}} \ d_{r-i,j}
\ S_{i,r+j,2r-i-j+1}\,,
\end{equation}
where the coefficients
$$
d_{r-i,j}=\langle j\moins 1,r\moins i\moins j \rangle
$$
are given by Eq. (\ref{Seg}) (or (\ref{eb})).
\end{theorem}

\bigskip

\noindent
We have the following values of ${\cal T}_1,
{\cal T}_2=\Phi({\cal T}_1)+{\overline{\cal T}_2}, \ldots,
{\cal T}_7=\Phi({\cal T}_6)+{\overline{\cal T}_7}$:
$$
\aligned
&S_{22}\\
&S_{133}\plus 3S_{34}\\
&S_{244}\plus 3S_{145}\plus 7S_{46}\plus 3S_{55}\\
&S_{355}\plus 3S_{256}\plus 7S_{157}\plus 3S_{166}\plus 15S_{58}
\plus 10S_{67}\\
&S_{466}\plus 3S_{367}\plus 7S_{268}\plus 3S_{277}\plus 15S_{169}\plus
10S_{178}\plus 31S_{6,10}\plus 25S_{79}\plus 10S_{88}\\
&S_{577}\plus 3S_{478}\plus 7S_{379}\plus 3S_{388}\plus 15S_{2,7,10}
\plus 10S_{289}\plus 31S_{1,7,11}\plus 25S_{1,8,10}\plus 10S_{199}\plus \\
&63S_{7,12}\plus 56S_{8,11}\plus 35S_{9,10}\\
&S_{688}\plus 3S_{589}\plus 7S_{4,8,10}\plus 3S_{499}\plus 15S_{3,8,11}
\plus 10S_{3,10,10}\plus 31S_{2,8,12}\plus 25S_{2,9,11}\plus\\
&10S_{2,9,10}\plus 63S_{1,8,13}\plus 56S_{1,9,12}\plus 35S_{1,10,11}
\plus 127S_{8,14}\plus 119S_{9,13}\plus 91S_{10,12}\plus 35S_{11,11}\,.
\endaligned
$$

\bigskip

\noindent
{\bf Acknowledgments}
Though the author of the present paper is responsible for the exposition
of the details, many computations here were done together with Alain
Lascoux. This help is gratefully acknowledged.
We also thank Richard Rimanyi for introducing us to his paper \cite{Rim2},
Andrzej Weber for his assistance on singularities, and Ozer Ozturk for
pointing out several defects of the manuscript.

\bigskip

\noindent
{\bf Note}
After completion of the first version \cite{P23} of this paper, we received
the preprint \cite{FR1} containing a Chern monomial expression for
the Thom series of $I_{2,2}$ supported by a computer evidence.
In July 2006, we have received a letter from Kazarian \cite{Ka2} informing
us that he has found another derivation of the Thom polynomial for $I_{2,2}$
based on geometric considerations.


\begin{thebibliography}{99}\small
\addcontentsline{toc}{section}{\string\numberline{}References}

\bibitem{AVGL} V. Arnold, V. Vasilev, V. Goryunov, O. Lyashko,
\emph{Singularities. Local and global theory,}
Enc. Math. Sci. vol. 6 (Dynamical Systems VI), Springer, 1993.

\bibitem{BFR} G. Berczi, L. Feher, R. Rimanyi,
\emph{Expressions for resultants coming from the global theory of
singularities,} in: ``Topics in algebraic and noncommutative geometry'',
(L. McEwan et al. eds.),
Contemporary Math. AMS {\bf 324} (2003), 63--69.

\bibitem{BR} A. Berele, A. Regev,
\emph{Hook Young diagrams with applications to combinatorics and to
representation theory of Lie superalgebras,}
Adv. in Math. {\bf 64} (1987), 118--175.

\bibitem{D} J. Damon,
\emph{Thom polynomials for contact singularities,}
Ph.D. Thesis, Harvard, 1972.

\bibitem{dPW} A. Du Plessis, C.T.C. Wall,
\emph{The geometry of topological stability,}
Oxford Math. Monographs, 1995.

\bibitem{FK} L. Feher, B. Komuves,
\emph{On second order Thom-Boardman singularities,}
Fund. Math. {\bf 191} (2006), 249--264.

\bibitem{FR} L. Feher, R. Rimanyi,
\emph{Calculation of Thom polynomials and other cohomological obstructions
for group actions,} in: ``Real and complex singularities (San Carlos
2002)'' (T. Gaffney and M. Ruas eds.), Contemporary Math. {\bf 354},
(2004), 69--93.

\bibitem{FR1} L. Feher, R. Rimanyi,
\emph{On the structure of Thom polynomials of singularities,}
Preprint (September 2005).

\bibitem{FP} W. Fulton, P. Pragacz,
\emph{Schubert varieties and degeneracy loci,}
Springer LNM {\bf 1689} (1998).

\bibitem{J} K. J\"anich,
\emph{Symmetry properties of singularities of $C^{\infty}$-functions,}
Math. Ann. {\bf 238} (1979), 147--156.

\bibitem{Ka} M.E. Kazarian,
\emph{Characteristic classes of singularity theory,}
in: ``The Arnold-Gelfand mathematical seminars: Geometry and singularity
theory'' (1997), 325--340.

\bibitem{Ka2} M. E. Kazarian,
\emph{Classifying spaces of singularities and Thom polynomials,}
in: ``New developments in singularity theory'', NATO Sci. Ser. II Math.
Phys. Chem., {\bf 21}, Kluwer Acad. Publ., Dordrecht (2001), 117--134.

\bibitem{Ka1} M. E. Kazarian,
\emph{Letter to the author, dated July 4, 2006.}

\bibitem{Kl} S. Kleiman,
\emph{The enumerative theory of singularities,}
in: ``Real and complex singularities, Oslo 1976'' (P. Holm ed.) (1978),
297--396.

\bibitem{LLT} D. Laksov, A. Lascoux, A. Thorup,
\emph{On Giambelli's theorem for complete correlations,}
Acta Math. {\bf 162} (1989), 143--199.

\bibitem{L} A. Lascoux,
\emph{Symmetric functions and combinatorial operators on polynomials},
CBMS/AMS Lectures Notes {\bf 99}, Providence (2003).

\bibitem{L1} A. Lascoux,
\emph{Addition of $\pm 1$: application to arithmetic},
S\'eminaire Lotharingien de Combinatoire, {\bf B52a} (2004), 9 pp.

\bibitem{LS} A. Lascoux, M-P. Sch\"utzenberger,
\emph{Formulaire raisonn\'e de fonctions sym\'e\-triques,}
Universit\'e Paris 7, 1985.

\bibitem{M} I.G. Macdonald,
\emph{Symmetric functions and Hall-Littlewood polynomials,}
Oxford Math. Monographs, Second Edition, 1995.

\bibitem{O} O. Ozturk,
\emph{On Thom polynomials for $A_4(-)$ via Schur functions,}
Preprint, IMPAN Warszawa 2006 (670) -- to appear in Serdica Math. J.

\bibitem{Po} I. Porteous,
\emph{Simple singularities of maps,}
in: ``Proc. Liverpool Singularities I'', Springer LNM {\bf 192} (1971),
286--307.

\bibitem{P1} P. Pragacz,
\emph{Note on elimination theory,}
Indagationes Math. {\bf 49} (1987), 215--221.

\bibitem{P} P. Pragacz,
\emph{Enumerative geometry of degeneracy loci,}
Ann. Sc. Ec. Norm. Sup. {\bf 21} (1988), 413--454.

\bibitem{P2} P. Pragacz,
\emph{Algebro-geometric applications of Schur $S$- and $Q$-polyno\-mials,}
in: ``Topics in invariant theory'' -- S\'eminaire d'Alg\`ebre
Dubreil-Malliavin 1989-1990 (M-P. Malliavin ed.), Springer LNM {\bf 1478}
(1991), 130--191.

\bibitem{P223} P. Pragacz
\emph{Symmetric polynomials and divided differences in formulas
of intersection theory,}
in: ``Parameter spaces'' (P. Pragacz ed.), Banach Center Publications
{\bf 36} (1996), 125--177.

\bibitem{P23} P. Pragacz,
\emph{Thom polynomials and Schur functions I,}
math.AG/0509234.

\bibitem{P3} P. Pragacz,
\emph{Thom polynomials and Schur functions: towards the singularities
$A_i(-)$,}
Preprint MPIM Bonn 2006 (139).

\bibitem{P4} P. Pragacz,
\emph{Thom polynomials and Schur functions: the singularities $A_3(-)$,}
in preparation.

\bibitem{PT} P. Pragacz, A. Thorup,
\emph{On a Jacobi-Trudi identity for supersymmetric polynomials,}
Adv. in Math. {\bf 95} (1992), 8--17.

\bibitem{PW} P. Pragacz, A. Weber,
\emph{Positivity of Schur function expansions of Thom polynomials,}
Preprint, math.AG/0605308 and MPIM Bonn 2006 (60) -- to appear in Fund. Math.

\bibitem{Rim2} R. Rimanyi,
\emph{Thom polynomials, symmetries and incidences of singularities,}
Inv. Math. {\bf 143} (2001), 499--521.

\bibitem{RS} R. Rimanyi, A. Sz\"ucs,
\emph{Generalized Pontrjagin-Thom construction for maps with
singularities,}
Topology {\bf 37} (1998), 1177--1191.

\bibitem{S} H. Schubert,
\emph{Allgemeine Anzahlfunctionen f\"ur Kegelschnitte, Fl\"achen und
Ra\"ume zweiten Grades in $n$ Dimensionen,}
Math. Ann. {\bf 45} (1894), 153--206.

\bibitem{Ste} J. Stembridge,
{A characterization of supersymmetric polynomials,}
J. of Algebra {\bf 95} (1985), 439--444.

\bibitem{T} R. Thom,
\emph{Les singularit\'es des applications diff\'erentiables,}
Ann. Inst. Fourier {\bf 6} (1955--56), 43--87.

\bibitem{W} C.T.C. Wall,
\emph{A second note on symmetry of singularities,}
Bull. London Math. Soc. {\bf 12} (1980), 347--354.

\end{thebibliography}
\end{document}